\documentclass{amsart}
\usepackage{hyperref}
\usepackage{amsthm}
\usepackage{amssymb}
\usepackage{amsmath}
\usepackage{amscd}
\usepackage{cases}
\usepackage[dvips]{graphicx}
\usepackage{latexsym}
\usepackage{pifont}

\newtheorem{thm}{Theorem}[section]

\newtheorem{prop}[thm]{Proposition}
\newtheorem{lem}[thm]{Lemma}

\newtheorem{pr}[thm]{Proof}

\makeatletter
 
 \@addtoreset{equation}{section}
\makeatother

\begin{document}

\title{{On positive and almost alternating links}}
\author{Kazuhiko Inoue}
\address{Graduate School of Mathematics, Kyusyu University, 744, Motooka, Nishi-ku, Fukuoka, 819-0395, Japan}

\maketitle

\begin{abstract}
In this paper, we show that a link which has a positive and almost alternating diagram is alternating, besides that a positive and non-alternating Montesinos link has an almost positive-alternating diagram.
\end{abstract}

\section{Introduction}
A \textbf{link} is a disjoint union of circles embedded in $\mathbb{S}^3$, and a \textbf{knot} is consist of one circle. A \textbf{diagram} of a link is a generic projection of a link on $\mathbb{S}^2$ with over/under information for each double point. A diagram is \textbf{alternating} if the over-crossings and under-crossings appear alternately along every component of the diagram, and a link is alternating if it has an alternating diagram. A link diagram is \textbf{almost alternating} if one crossing change makes it into an alternating diagram, and a link is almost alternating if it has an almost alternating diagram and no alternating diagram.

A diagram is \textbf{positive} if the sign of every crossing is positive. (A negative diagram is the mirror image of a positive diagram.) A link is positive if it has a positive diagram, and a link is \textbf{positive and alternating} if it has a positive diagram and an alternating diagram. A link is \textbf{positive-alternating} if it has a positive and alternating diagram. Nakamura showed that every positive and alternating link has a positive-alternating diagram (\cite{Na}). (We call a positive-alternating diagram, PA-diagram, and a positeve and alternating link, PA-link for short.) So our concern is \textbf{a positive and almost alternating link}, that is to say a link which has a positive diagram and almost alternating diagram and has no alternating diagram. In section 3, we show the following:

\textbf{Theorem\ref{thm:111}} Let $L$ an oriented link. If $L$ has a positive and almost alternating diagram then $L$ is alternating.

Besides we know that every positive and almost alternating knot has an almost positive-alternating diagram with up to eleven crossings. Furthermore Jong and Kishimoto showed that every positive knot up to genus two is  positive-alternating or almost positive-alternating(\cite{Jo}). A diagram is \textbf{almost positive-alternating} if one crossing change makes it into a PA-diagram. We say such diagram almost PA-diagram. In this paper we show the following:

\textbf{Proposition\ref{prop:112}} Every positive Montesinos link has an almost positive-alternating diagram.

This paper is organized as follows; In \textbf{Section 2}, we briefly review almost alternating links and positive alternating links. In \textbf{Section 3}, we prove Theorem\ref{thm:111}. In \textbf{Section 4}, we characterize positive Montesinos links and prove Proposition\ref{prop:112}.

\section{Preliminary}

\begin{figure}[htbp]
  \centering
  \includegraphics[width=9cm,clip]{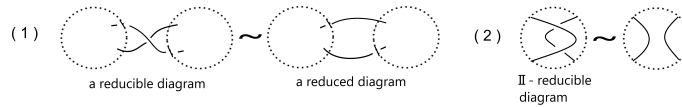}
  \caption{A reducible diagram and a II-reducible diagram}
  \label{fig:2-2ji}
 \end{figure}

First we shortly introduce some definitions. A diagram $D$ is said to be \textbf{equivalent} to a diagram $D'$ if they both represent same link. A diagram is said to be \textbf{reduced} if there exists no crossing such that the diagram is separated by splicing the crossing as shown in Figure\ref{fig:2-2ji}(1). A diagram is said to be \textbf{II-reduced} if there are no obvious removal Reidemeister-II move i.e. the link contains no 2-tangle as shown in Figure \ref{fig:2-2ji}(2) (See \cite{Ts}).

\begin{figure}[htbp]
  \centering
  \includegraphics[width=5cm,clip]{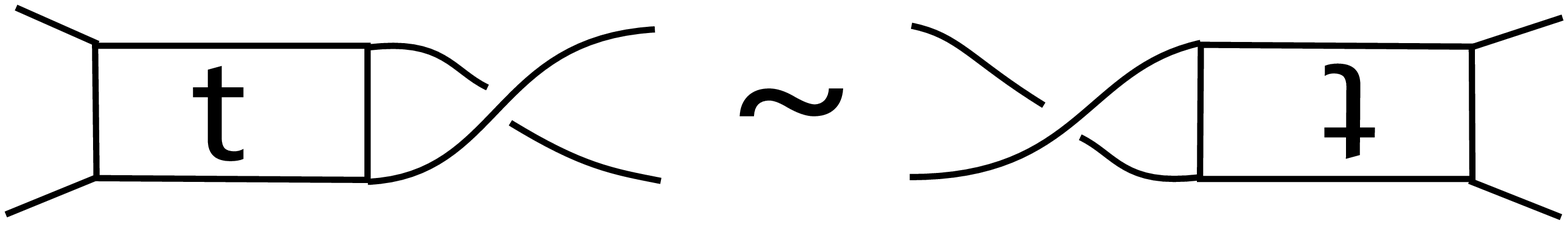}
  \caption{A flype}
  \label{fig:2-3ji}
 \end{figure}

A \textbf{flype} is an isotopy move applied on a sub tangle of the form $[ \pm1 ] + t$, and it fixes the endpoints of the sub tangle. See Figure \ref{fig:2-3ji}. A flype preserves the alternating structure of a diagram.(\cite{Ka})

\begin{figure}[htbp]
  \centering
  \includegraphics[width=4cm,clip]{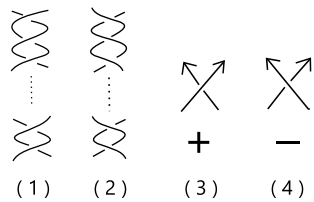}
  \caption{A positive(negative) tangle and the sign of a crossing point}
  \label{fig:2-4ji}
 \end{figure}

We distinguish a positive tangle from a positive diagram as the following. A tangle is positive if it is as shown in   Figure \ref{fig:2-4ji}(1) and negative as shown in (2). The sign of a crossing point is $+$ if it is as shown in Figure \ref{fig:2-4ji}(3) and $-$ as shown in (4). A diagram is positive (resp. negative) if every crossing point in the diagram has the same sign $+$ (resp. $-$).

\begin{figure}[htbp]
  \centering
  \includegraphics[width=9cm,clip]{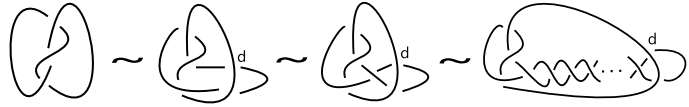}
  \caption{Almost alternating diagrams of a trefoil knot}
  \label{fig:2-1ji}
 \end{figure}

Next we introduce some results about almost alternating links. First any alternating link involving a trivial link has an almost alternating diagram. Moreover any alternating link has infinite almost alternating diagrams(see \cite{Ad}). For example we can make infinite almost alternating diagrams from a trefoil knot as shown in Figure \ref{fig:2-1ji}. Every diagram turns into an alternating diagram if we change the crossing point $d$.

\begin{thm}\label{thm:113}
 
Every positive and almost alternating knot is almost positive-alternating with up to eleven crossings $($\cite{Cr}$)$.

\end{thm}
 
\begin{thm}\label{thm:114}

Positive knots up to genus two are positive-alternating or almost positive-alternating $($\cite{Jo}$)$.

\end{thm}

\begin{thm}\label{thm:115}

Non-alternating Montesinos links are almost alternating $($\cite{Ab}$)$.

\end{thm}

\begin{thm}\label{thm:116}

Any reduced alternating diagram of a positive alternating link is positive-alternating$($\cite{Na}$)$.

\end{thm}

Then our concern at the moment is the following question.

Question: How is the diagram of a positive and almost alternating link? 

 In section 3, we show that a positive and almost alternating link does not have a positive and almost alternating diagram for the partial answer to the question above.

\section{Main Theorem}

\begin{thm}\label{thm:111}

 Let $L$ be an oriented link. If $L$ has a positive and almost alternating diagram then $L$ is alternating.

\end{thm}

\begin{pr}\label{pr117}
 By the assumption above $L$ has an almost alternating diagram, so we see $L$ is alternating or almost alternating. Our  
 claim is that every positive and almost alternating diagram of $L$ is equivalent to an alternating diagram. First of all if a diagram $D$ is reducible then $D$ is equivalent to an alternating diagram as shown in Figure \ref{fig:3-1ji}. If $D$ is {\rm II}-reducible we can see this is equivalent to an alternating diagram in a similar fashion. Therefore we can assume that $D$ is both reduced and {\rm II}-reduced.

\begin{figure}[htbp]
  \centering
  \includegraphics[width=5cm,clip]{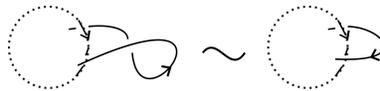}
  \caption{A reducible almost alternating diagram}
  \label{fig:3-1ji}
\end{figure}

 In general positive and almost alternating diagrams are as shown in Figure\ref{fig:3-2ji}, where the diagram in every shaded portion is positive-alternating.

\begin{figure}[htbp]
  \centering
  \includegraphics[width=12cm,clip]{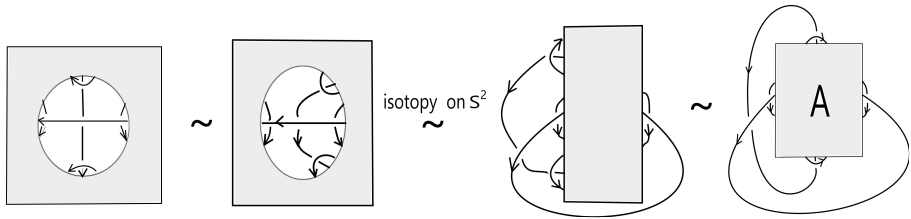}
  \caption{positive and almost alternating diagram}
  \label{fig:3-2ji}
\end{figure}

The shaded portion in the rightmost figure is equivalent to a disk and we denote this region by $A$. Note that there does not happen the case as shown in Figure\ref{fig:3-3ji}, because this diagram is not positive. 

\begin{figure}[htbp]
  \centering
  \includegraphics[width=3.cm,clip]{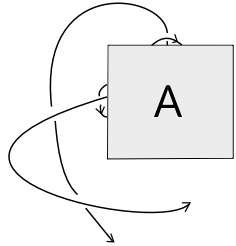}
  \caption{This diagram is not positive}
  \label{fig:3-3ji}
\end{figure}

Assume disk $A$ separates into disk $A_1$ and disk $A_2$. Since each diagram in $A_1$ and $A_2$ is alternating, then  the diagram $D$ is equivalent to an alternating diagram $D'$. $(See Figure \ref{fig:3-4ji}.$) Hence $L$ is alternating. 

\begin{figure}[htbp]
  \centering
  \includegraphics[width=12cm,clip]{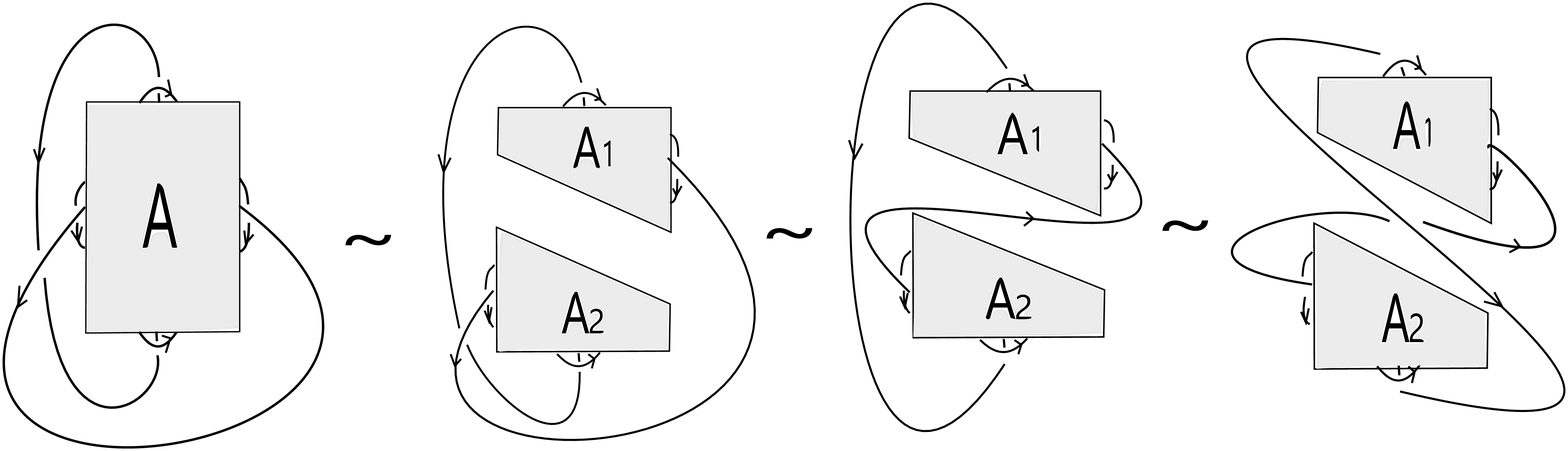}
  \caption{}
  \label{fig:3-4ji}
\end{figure}

Next we prove that disk $A$ actually separates into disk $A_1$ and disk $A_2$. We name five crossing points outside of $A$, $\alpha$, $\alpha'$, $\beta$, $\beta'$, $d$ as shown in Figure  \ref{fig:3-4-2ji}. Besides we also name the strand which passes through $\alpha$ and enters into $A$, $\overline{\alpha}$, similarly the strand which passes through $\beta$ and enters into $A$, $\overline{\beta}$.

\begin{figure}[htbp]
  \centering
  \includegraphics[width=3cm,clip]{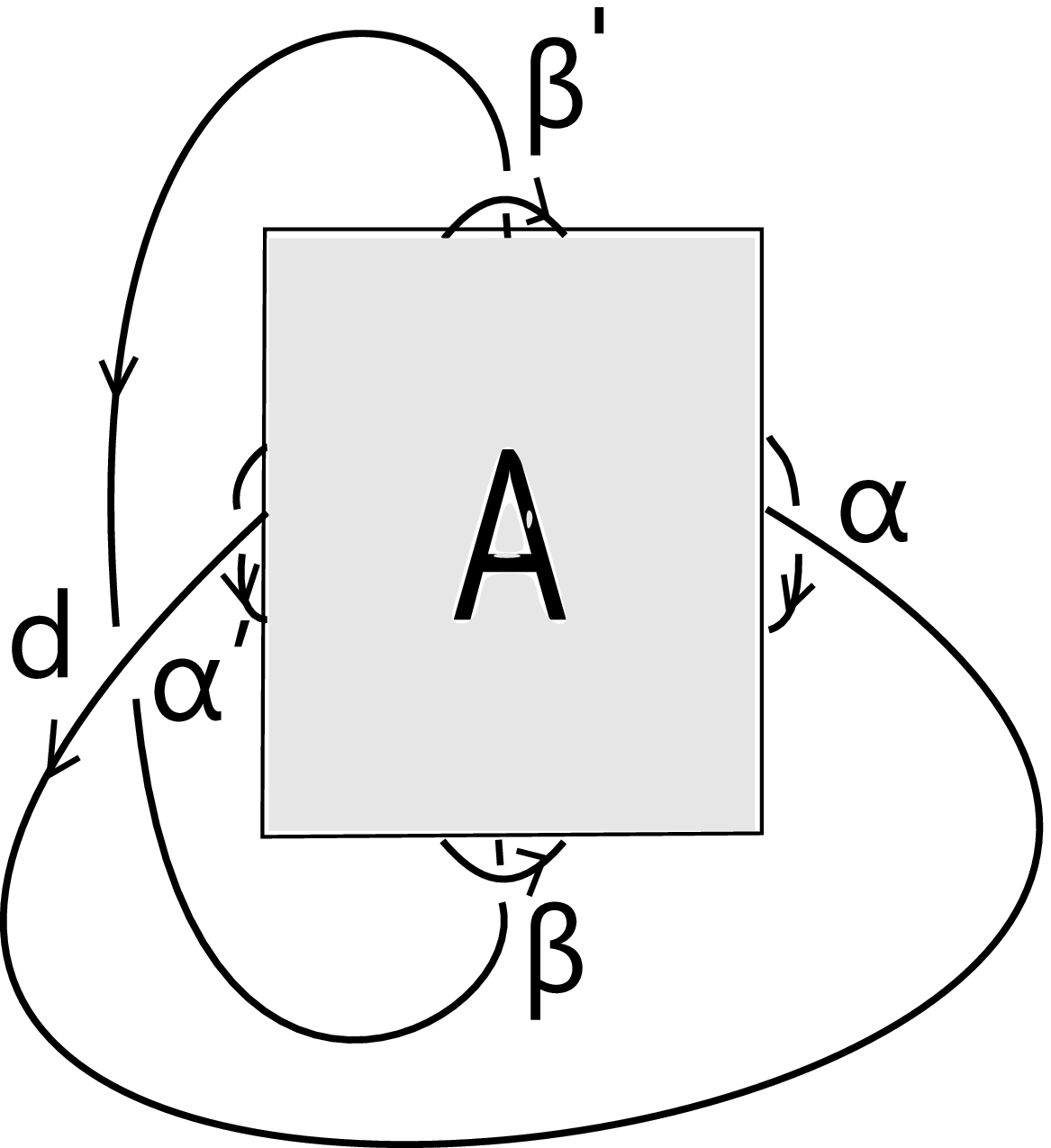}
  \caption{}
  \label{fig:3-4-2ji}
\end{figure}

When the strand which passed under the strand $\overline{\alpha}$ at $\alpha$ crosses next strand, there can be three cases as  shown in Figure \ref{fig:3-5ji} $(1) \sim (3)$. 

\begin{figure}[htbp]
  \centering
  \includegraphics[width=5cm,clip]{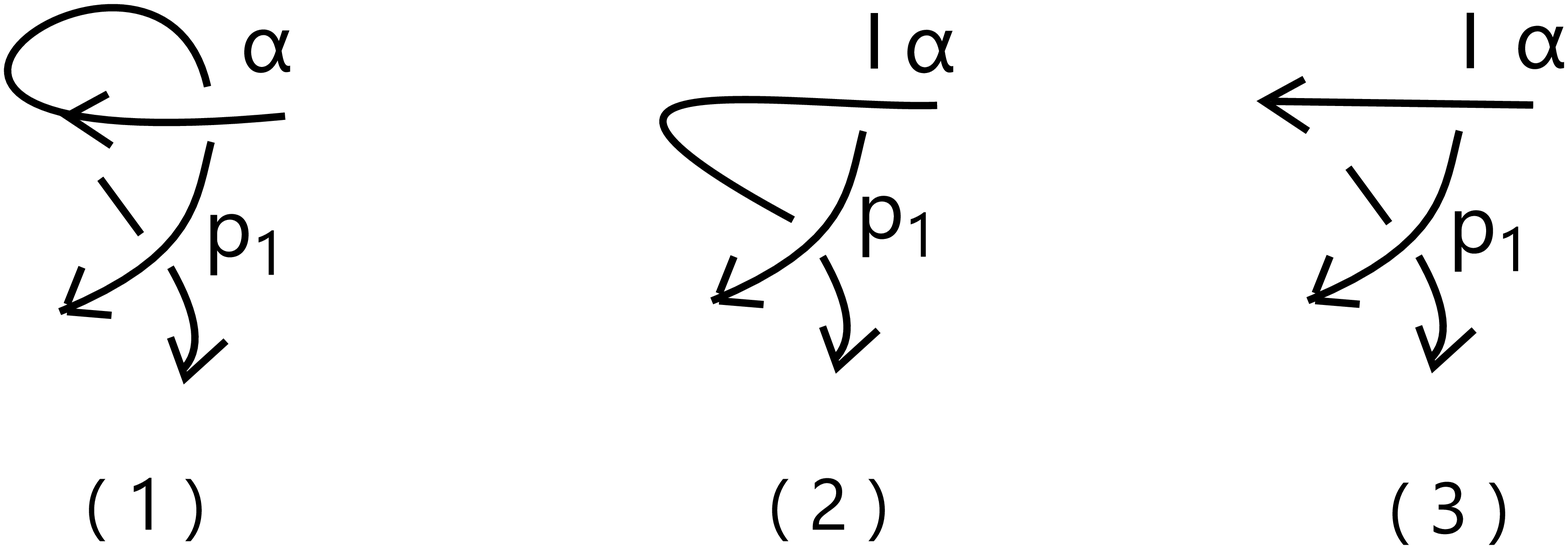}
  \caption{}
  \label{fig:3-5ji}
\end{figure}

 Since the diagram $D$ is positive and alternating in region $A$ therefore in any case the next strand passes under this strand from the right side to the left side as shown in Figure \ref{fig:3-6ji}$(1)$. We name these crossing points $p_1, p_2, \ldots$ and also name the arc from $\alpha$ to $p_1$, $\overline{p}_0$, from $p_1$ to $p_2$, $\overline{p}_1$, similarly $\overline{p}_2$,  $\overline{p}_3$, $\ldots$ and so on. On the other hand we consider the strand which passes over  $\overline{\beta}$ as shown in Figure \ref{fig:3-6ji}$(2)$.

\begin{figure}[htbp]
  \centering
  \includegraphics[width=8cm,clip]{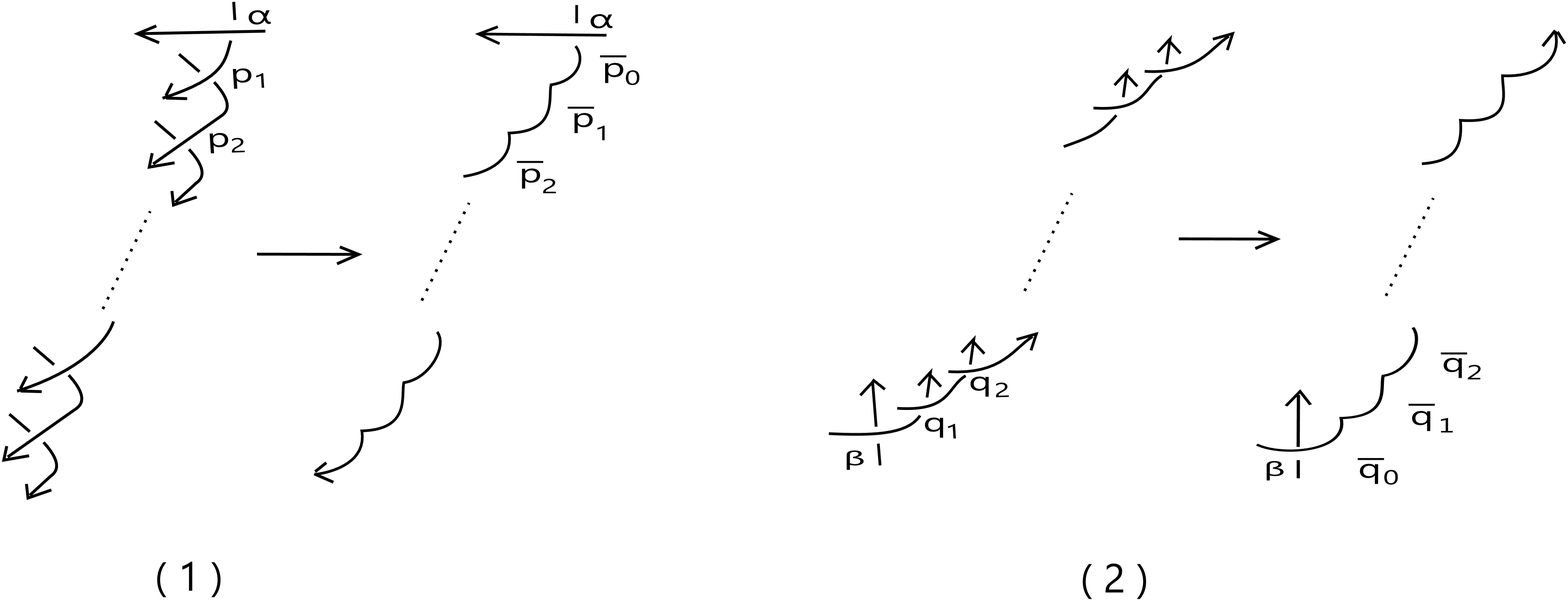}
  \caption{The strand which passes under $\overline{\alpha}$ and the strand which passes over $\overline{\beta}$}
  \label{fig:3-6ji}
\end{figure}

In the case where a strand crosses a loop or a strand crosses by itself, we regard as shown in Figure \ref{fig:3-6-2ji}.

\begin{figure}[htbp]
  \centering
  \includegraphics[width=6cm,clip]{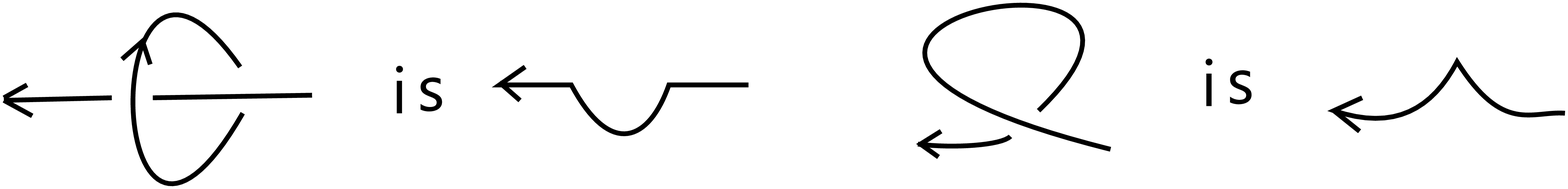}
  \caption{}
  \label{fig:3-6-2ji}
\end{figure}

Finally there are two sequences of arcs in $A$ and they are both oriented. This is such as shown in Figure \ref{fig:3-7ji}$(1)$. If $\overline{p}_m$ and $\overline{q}_n$ cross each other then  $\overline{p}_m$ passes over $\overline{q}_n$ from the left side to the right side as shown in Figure \ref{fig:3-7ji}$(2)$.

\begin{figure}[htbp]
  \centering
  \includegraphics[width=9cm,clip]{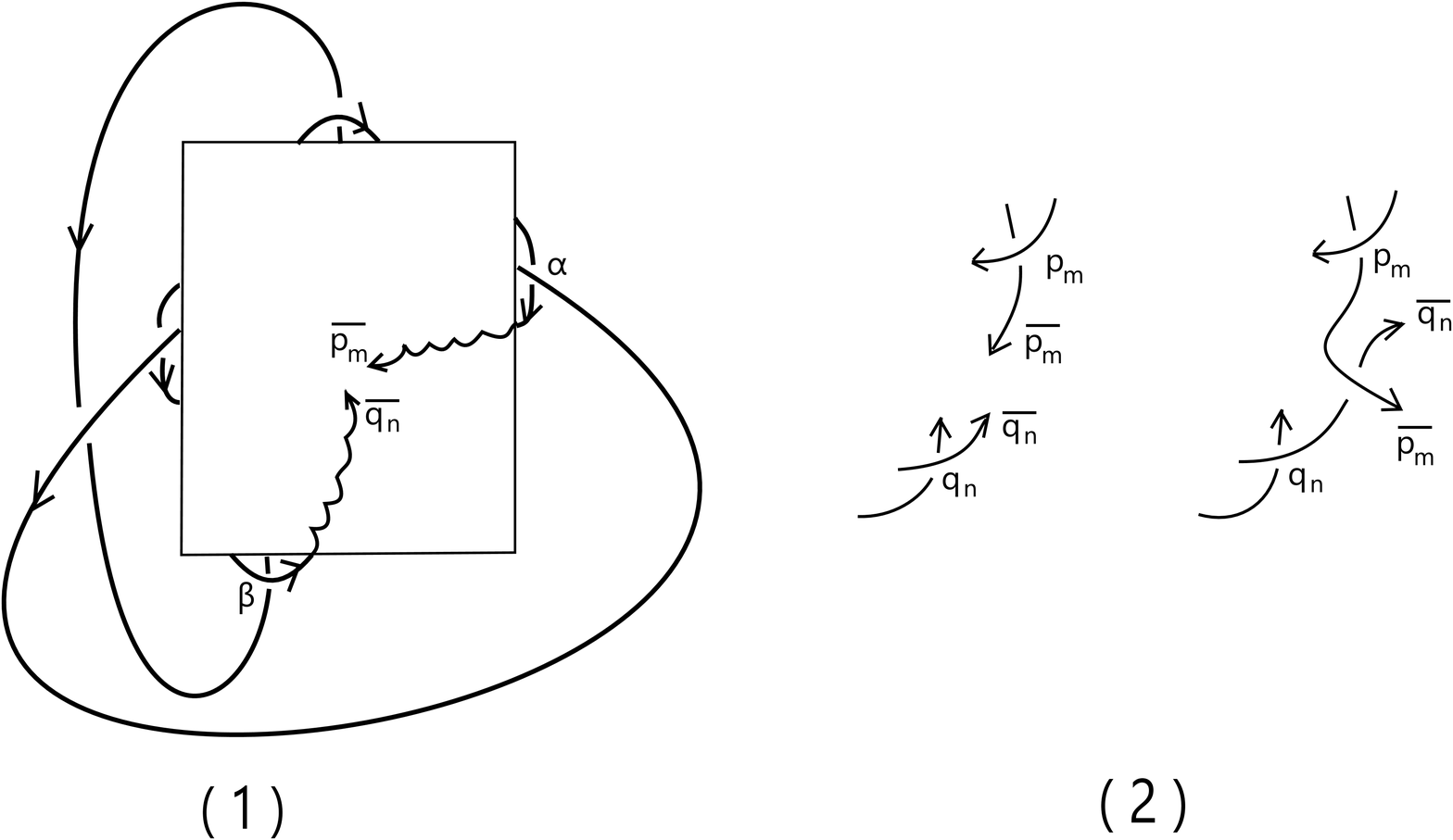}
  \caption{The relation between $\overline{p}_m$ and $\overline{q}_n$ }
  \label{fig:3-7ji}
\end{figure}

We name this crossing point $c$, then there is a polygon with vertices $\alpha$, $d$, $\beta$, $q_1$, $q_2$, $\ldots$, $q_n$, $c$, $p_m$, $p_{m-1}$, $\ldots$, $p_2$, $p_1$. And two arcs $\overline{p}_m$, $\overline{q}_n$ enter this polygon as     shown in Figure\ref{fig:3-9ji}. This is the contradiction to \textbf{The Jordan curve theorem}$($\cite{Ha}$)$. 

\begin{thm}\label{thm:118}$($Jordan curve theorem$)$

Let $C$ be the image of the unit circle, that is $C = \{ ( x, y ); x^2 + y^2 = 1 \}$ under an injective continuous mapping  $\gamma$ into $\mathbb{R}^2$. Then $\mathbb{R}^2 \backslash C$ is disconnected and consists of two component.

\end{thm}

 Moreover if $\overline{p}_m$ or $\overline{q}_n$ crosses some arc in $\{\overline{p}_i\}$ or $\{\overline{q}_j\}$ then next it crosses the same arc and enter this polygon again. Because each $\overline{p}_i$ is an arc from under crossing to over crossing and each $\overline{q}_j$ is an arc from over crossing to under crossing. After all we can see that $\overline{p}_m$ and $\overline{q}_n$ never cross each other in $A$.

\begin{figure}[htbp]
  \centering
  \includegraphics[width=5cm,clip]{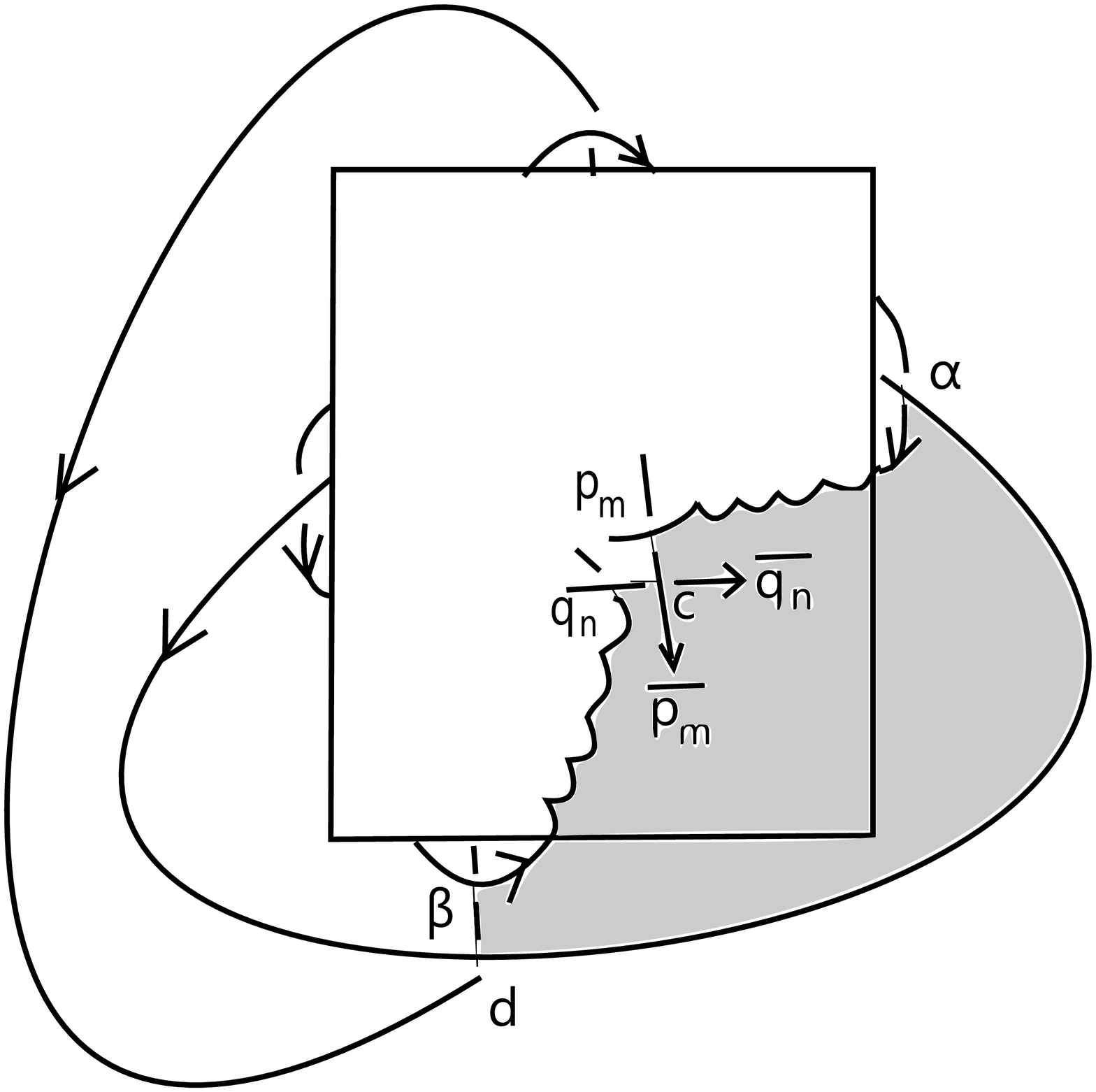}
  \caption{}
  \label{fig:3-9ji}
\end{figure}

For this reason $A$ must separate into $A_1$ and $A_2$ hence $D$ is equivalent to an alternating diagram $D'$. This completes the proof of Theorem3.1.

\qed\end{pr}

From the theorem above we know that a positive and almost alternating link has no positive and almost alternating diagram. Furthermore we think the question again. How is the diagram of a positive and almost alternating link? And we give a partial result of this question in section 4.

\section{Positive and almost alternating Montesinos link}

\begin{figure}[htbp]
  \centering
  \includegraphics[width=11cm,clip]{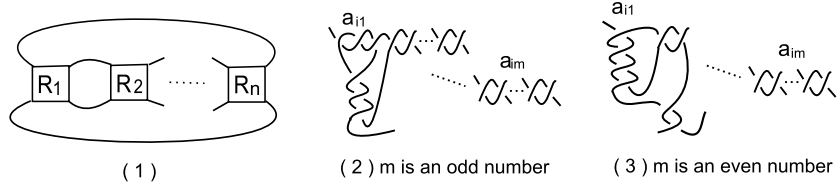}
  \caption{The standard diagram of a Montesinos link}
  \label{fig:4-1-1ji}
\end{figure}

In this section we would like to study an oriented Montesinos link $L$ denoted by
 $C( {{\alpha}_1}/{{\beta}_1}, {{\alpha}_2}/{{\beta}_2}, \ldots, {{\alpha}_n}/{{\beta}_n} )$, ${{\alpha}_i}/{{\beta}_i} \in \mathbb{Q}$. Any ${{\alpha}_i}/{{\beta}_i}$ represents not only a rational number but also a rational tangle $R_i = ({{\alpha}_i}/{{\beta}_i})$. About rational tangles, see \cite{Mu}. The standard diagram $D$ of $L$ denoted by $D({{\alpha}_1}/{{\beta}_1}, {{\alpha}_2}/{{\beta}_2}, \ldots, {{\alpha}_n}/{{\beta}_n})$ is shown in Figure\ref{fig:4-1-1ji}, where $({{\alpha}_i}/{{\beta}_i}) = R_i = R(a_{i1}, a_{i2}, \ldots, a_{im})$. That is to say, $D$ is the numerator of the sum of n rational tangles. For example in the case where any $a_{ij} > 0$ $R_i$ is as shown in Figure \ref{fig:4-1-1ji}(2) or (3). 

Abe and Kishimoto showed that any non-alternating Montesinos link is almost alternating, and we have the following proposition.

\begin{prop}\label{prop:112}

Let $L$ be an oriented Montesinos link and be denoted by 

$C({{\alpha}_1}/{{\beta}_1}, {{\alpha}_2}/{{\beta}_2}, \ldots, {{\alpha}_n}/{{\beta}_n})$ where ${{\alpha}_i}/{{\beta}_i} \in \mathbb{Q}$, and $D$ the standard diagram of $L$ such that $D({{\alpha}_1}/{{\beta}_1}, {{\alpha}_2}/{{\beta}_2}, \ldots, {{\alpha}_n}/{{\beta}_n})$. If $D$ is positive then $L$ has an almost PA-diagram.

\end{prop}

It is to be noted that in general if a link $L$ has a PA-diagram, then $L$ has also an almost PA-diagram. Because we can transform a PA-diagram $D$ of $L$ into an almost PA-diagram $D'$. See Figure \ref{fig:3-1ji}.

Before prooving the proposition above, we prove two other propositions and one lemma needed later. 

\begin{prop}\label{prop:113}

Let $L$ be an oriented link and $D$ be a diagram of $L$ such that $D = D_1 \sharp D_2 \sharp \cdots \sharp D_m$ where any $D_t$ is an alternating diagram $(1 \leq t \leq m)$. If $D$ is positive, then $L$ has a PA-diagram.

\end{prop}

\begin{figure}[htbp]
  \centering
  \includegraphics[width=6cm,clip]{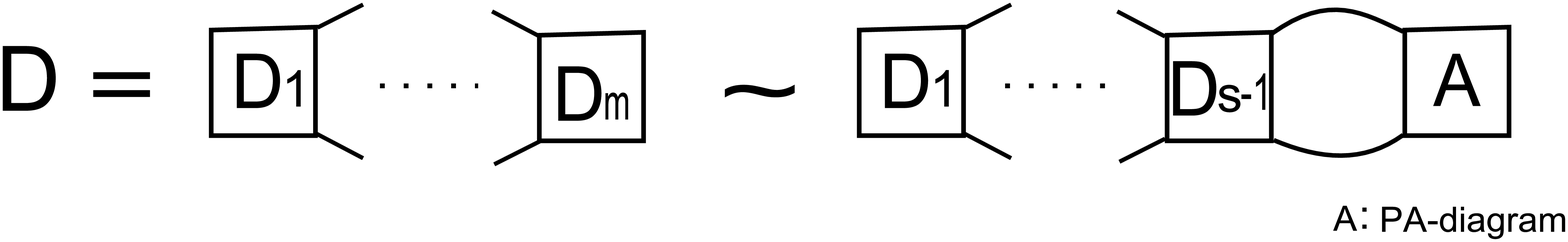}
  \caption{ }
  \label{fig:4-1-7ji}
\end{figure}

\begin{figure}[htbp]
  \centering
  \includegraphics[width=10cm,clip]{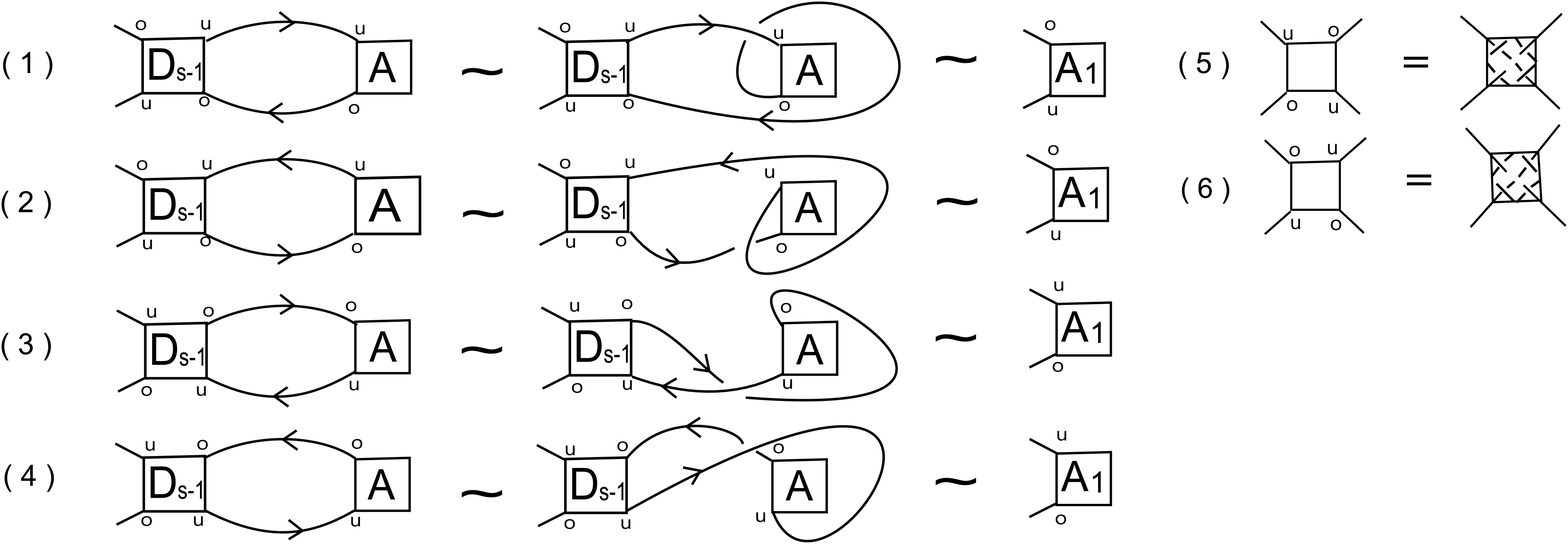}
  \caption{Four conditions of $D_{s-1} \sharp A$}
  \label{fig:4-1-3ji}
\end{figure}

\begin{figure}[htbp]
  \centering
  \includegraphics[width=6cm,clip]{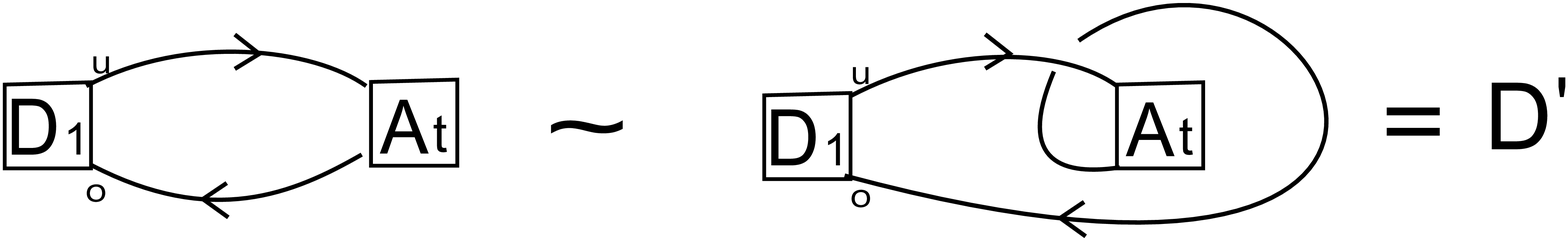}
  \caption{ }
  \label{fig:4-1-15ji}
\end{figure}

\begin{pr}\label{pr:120}

Assume $D = D_1 \sharp D_2 \sharp \cdots \sharp D_m$ is positive and $A =  D_s \sharp D_{s+1} \sharp \cdots \sharp D_m$ $(1 \leq s \leq m)$ is alternating as shown in Figure \ref{fig:4-1-7ji}. We consider the directions of two arcs on the left-hand side of $A$ and the over/under informations of the leftmost crossings of $A$ and the rightmost crossings of $D_{s-1}$. Then we can see the four conditions as shown in Figure \ref{fig:4-1-3ji}, where the simbol $o$ $($resp. $u$$)$ means that an over-crossing $($resp.  under-crossing$)$ appears first when we traverse the component from the end point $($\cite{Ab}$)$. By repeating this transformation we can finally obtain a PA-diagram of $L$ as shown in Figure \ref{fig:4-1-15ji}.

\qed\end{pr}

\begin{prop}\label{prop:114}

Let $L$ be an oriented Montesinos link, and $D$ be the standard diagram of $L$ denoted by $D( {{\alpha}_1}/{{\beta}_1}, {{\alpha}_2}/{{\beta}_2}, \ldots, {{\alpha}_n}/{{\beta}_n})$, \ ${{\alpha}_i}/{{\beta}_i} \in \mathbb{Q}$. Assume that $D$ is positive,  $\left| {{\alpha}_i}/{{\beta}_i} \right| \geq 1$ and ${\beta}_i \neq 0$ for any $i$ \ $(1 \leq i \leq n)$. Then $D$ is alternating.

\end{prop}

\begin{figure}[htbp]
  \centering
  \includegraphics[width=8cm,clip]{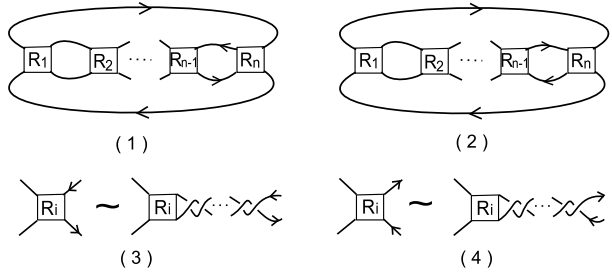}
  \caption{The case where the directions of lefthand side arcs of $R_1$ are oppsite}
  \label{fig:4-1-6ji}
\end{figure}

\begin{pr}\label{pr:121}

First we consider the case where the directions of left-hand side arcs of $R_1$ are parallel. In this case, naturally the directions of the right-hand side arcs of $R_n$ are also parallel as shown in Figure \ref{fig:4-1-5ji}$(1)$. Besides it is easy to see these directions hold in the case of $R_{n-1}, R_{n-2}, \ldots, R_2$ as shown in Figure \ref{fig:4-1-5ji}$(2)$. That is to say, the directions of the left-hand side arcs of any tangle $R_i$ are all the same as shown in Figure \ref{fig:4-1-5ji}$(3)$. Since each tangle $R_i = R(a_{i1}, a_{i2}, \ldots, a_{im})$ is positive and alternating, we know that $a_{ij} \leq 0$ for any $j$ for any $j$ $(1 \leq j \leq m)$. Hence ${{\alpha}_i}/{{\beta}_i} < 0$ for any $i$ $(1 \leq i \leq n)$. Then $D$ is necessarily alternating.

\begin{figure}[htbp]
  \centering
  \includegraphics[width=9cm,clip]{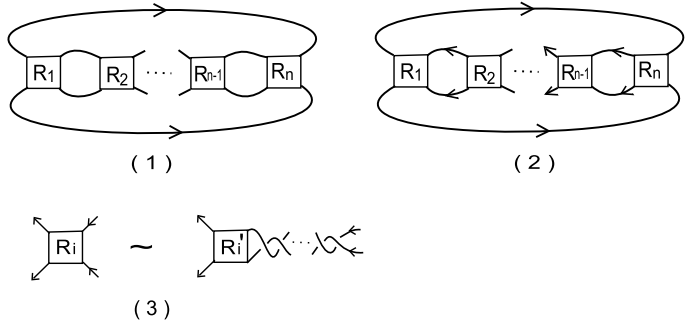}
  \caption{The case where the directions of lefthand side arcs of $R_1$ are parallel}
  \label{fig:4-1-5ji}
\end{figure}

Next we consider the case where the directions of the left-hand side arcs of $R_1$ are opposite. In this case, the directions of the right-hand side arcs of $R_n$ are as shown in Figure \ref{fig:4-1-6ji}$(1)$ or $(2)$. So for any tangle $R_i$, the directions of the right-hand side arcs are as shown in Figure \ref{fig:4-1-6ji}$(3)$ or $(4)$. In any case, we know that $a_{ij} > 0$ for any $j$ $(1 \leq j \leq m)$, because any $R_i$ is positive and alternating. Therefore ${{\alpha}_i}/{{\beta}_i} > 0$ for any $i$ $(1 \leq i \leq n)$ and $D$ must be positive. This completes the proof of the proposition.

\qed\end{pr}

\begin{figure}[htbp]
  \centering
  \includegraphics[width=7cm,clip]{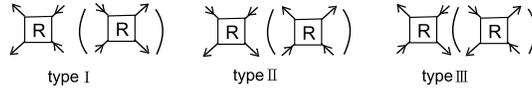}
  \caption{Three types of oriented tangles}
  \label{fig:4-1-8ji}
\end{figure}

In addition when we meditate upon oriented rational tangles, we can classify them into three types as shown in Figure \ref{fig:4-1-8ji}. What is more we can have the next lemma.

\begin{lem}\label{lem:122}

Let $R$ be an oriented rational tangle denoted by $({\alpha}/{\beta})$, where ${\alpha}/{\beta} \in \mathbb{Q}_{\neq 0}$, ${\beta} \neq 0$. If any crossing in $R$ has the same sign $+$, then the following holds.

\begin{enumerate}

\item If $R$ is of type {\rm I}, then ${\alpha}/{\beta} < 0$.

\item If $R$ is of type {\rm II}, then ${\alpha}/{\beta} > 0$.

\item If $R$ is of type {\rm III} and $\left| {\alpha}/{\beta} \right| \geq 1$, then ${\alpha}/{\beta} > 0$.

\item If $R$ is of type {\rm III} and $\left| {\alpha}/{\beta} \right| < 1$, then ${\alpha}/{\beta} < 0$.

\end{enumerate}
 
\end{lem}

\begin{figure}[htbp]
 \centering
 \includegraphics[width=11cm,clip]{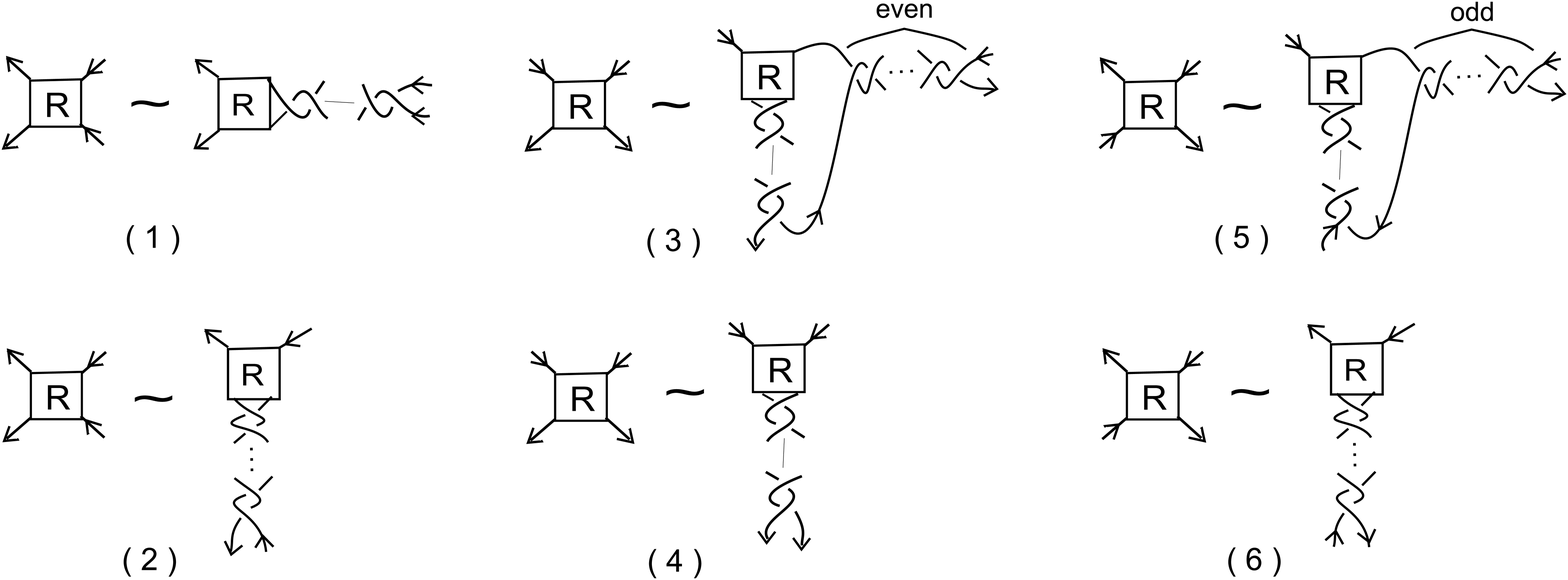}
 \caption{}
 \label{fig:4-1-11ji}
\end{figure}

\begin{pr}\label{pr:122}

In the case where $R$ is of type {\rm I}, the oriented tangle $R$ is naturally as shown in Figure \ref{fig:4-1-11ji}$(1)$ or $(2)$, and in both cases ${\alpha}/{\beta} < 0$. If we reverse all directions, we can prove in exactly the same way. In the case where $R$ is of type {\rm II}, $R$ is as shown in  Figure \ref{fig:4-1-11ji}$(3)$ or $(4)$, and it is easy to see in both cases  ${\alpha}/{\beta} > 0$. Besides, when $R$ is of type {\rm III} and $\left| {\alpha}/{\beta} \right| \geq 1$, $R$ is neccessarily as shown in Figure \ref{fig:4-1-11ji}$(5)$, and  ${\alpha}/{\beta} > 0$. On the contrary if $\left| {\alpha}/{\beta} \right| < 1$, $R$ must be as shown in Figure \ref{fig:4-1-11ji}$(6)$, and ${\alpha}/{\beta} < 0$. We have thus proved the lemma.

\qed\end{pr}

In fact, there are two types in type III as shown in  Figure \ref{fig:4-1-11ji}(5) and (6). So next we rename type III as shown in Figure \ref{fig:4-1-11ji}(5) type ${\rm III}_+$, and as shown in  Figure \ref{fig:4-1-11ji}(6) type ${\rm III}_-$.
Now we are ready to prove Proposition 4.1.

\begin{figure}[htbp]
  \centering
  \includegraphics[width=11cm,clip]{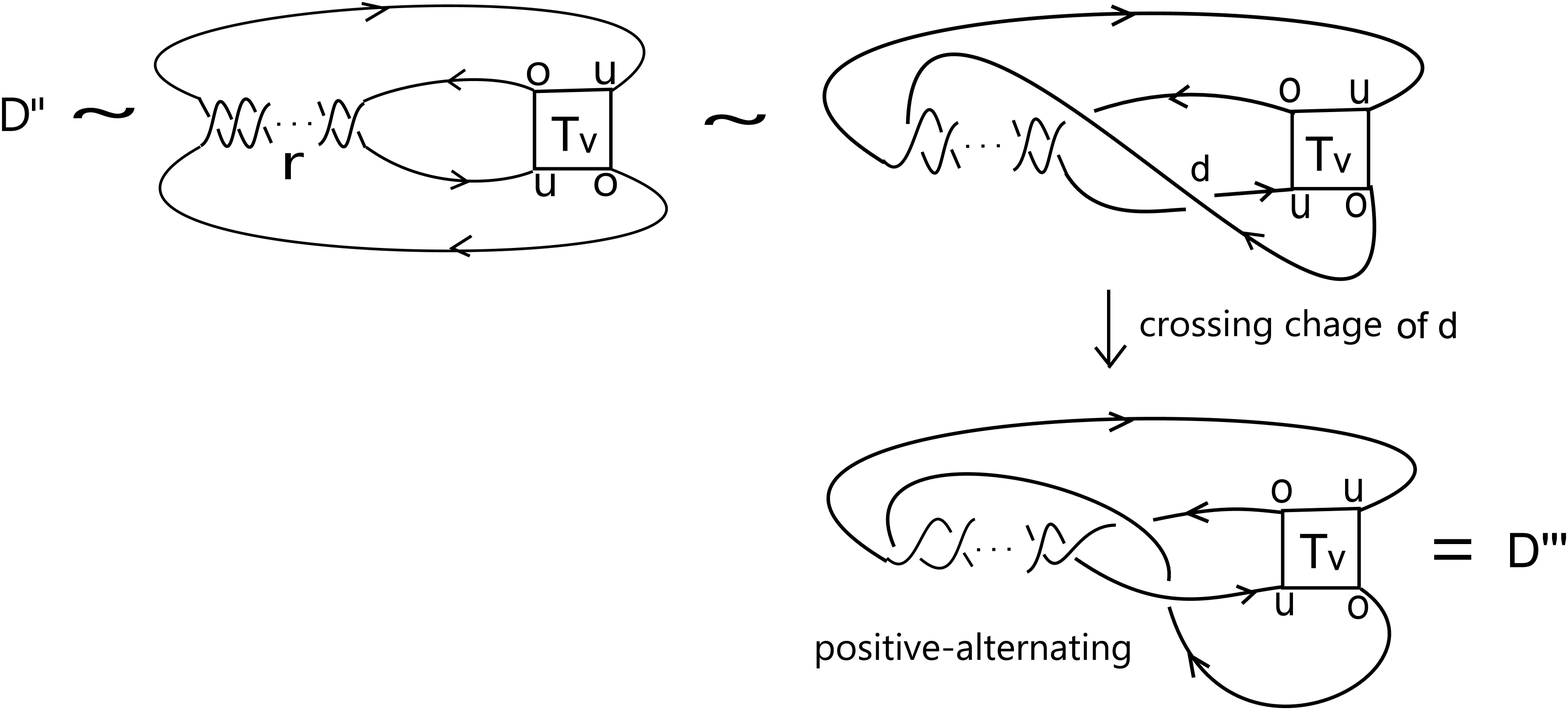}
  \caption{}
  \label{fig:4-1-14ji}
\end{figure}

\begin{pr}\label{pr:123} $(${\rm Proof of Proposition 4.1}$)$

First we consider the case where some ${\beta}_j = 0$. In this case the tangle ${{\alpha}_j}/{{\beta}_j}$ is a $\infty$-tangle as shown in Figure \ref{fig:4-1-4ji}$(1)$, and the diagram $D$ is like as shown in Figure \ref{fig:4-1-4ji}$(2)$, where each $R_k (1 \leq k \leq j-1, \ j+1 \leq k \leq m)$ is an alternating tangle. So by thinking that the denominator of a tangle $R_k$ is equivalent to a diagram $D_k$, we can regard $D = D_{j+1} \sharp \cdots \sharp D_m \sharp D_1 \sharp \cdots \sharp D_{j-1}$, where any $D_k$ is alternating. Therefore if $D$ is positive then $L$ has a PA-diagram by Proposition $4.2$. Hence $L$ has also an almost PA-diagram.

\begin{figure}[htbp]
  \centering
  \includegraphics[width=10cm,clip]{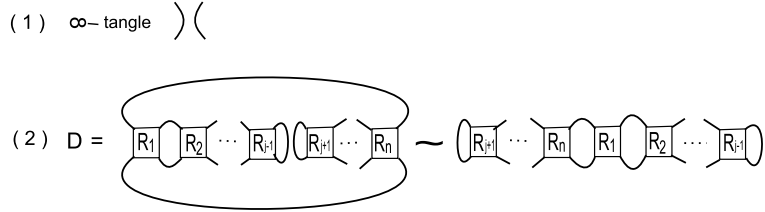}
  \caption{The case where some $\beta_j = 0$}
  \label{fig:4-1-4ji}
\end{figure}

Next we consider the case where $\beta_i \neq 0$ for any $i$ \ $(1 \leq i \leq n)$. In this case, from Proposition 4.2 if $\left| {\alpha_i}/{\beta_i} \right| \geq 1$ for any $i$ \ $(1 \leq i \leq n)$ then $D$ is alternating. Thus if $D$ is non-alternating then there must be some ${\alpha_j}/{\beta_j}$ such that  $\left| {{\alpha}_j}/{{\beta}_j} \right| < 1$. That is to say, there exists some tangle $R_j$ such that  $R_j = R(a_{j1}, a_{j2}, \ldots, a_{jm})$, $a_{jm} = 0$. Futhermore, by Lemma 4.6 we know that there exist some $($may be one$)$ rational tangles of type ${\rm III}_-$ and some $($may be one$)$ rational tangles of type ${\rm II}$ or ${\rm III}_+$ in $\{R_i\}$. In this condition we can transform the tangles of type {\rm II} as shown in Figure \ref{fig:4-1-11ji}$(3)$ into such tangles as shown in Figure \ref{fig:4-1-9ji}$(1)$ $($or all the directions are opposite$)$, and the tangles of type ${\rm III}_+$ as shown in Figure \ref{fig:4-1-11ji}$(5)$ into the tangles as shown in Figure \ref{fig:4-1-9ji}$(2)$ $($or all the directions are opposite$)$.

\begin{figure}[htbp]
 \centering
 \includegraphics[width=12cm,clip]{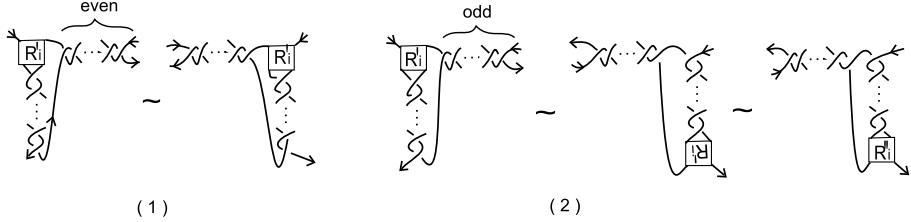}
 \caption{The transformation of oriented tangles}
 \label{fig:4-1-9ji}
\end{figure}

After the transformations above, we can regard that $D$ is equivalent to the diagram $D' = D(P_1, P_2, \ldots, P_n)$ as shown in Figure \ref{fig:4-1-10ji}$(1)$ where $P_k (1 \leq k \leq n)$ is a rational tangle or a 180 degree reversed rational tangle as shown in Figure \ref{fig:4-1-10ji}$(2) \sim (5)$. Besides, there exists at least one tangle of type ${\rm III}_-$, that is, the tangles as shown in Figure \ref{fig:4-1-10ji}$(4)$ or $(5)$. These tangles as shown in Figure \ref{fig:4-1-10ji}$(2) \sim (5)$ are all alternating and all crossings in these diagram have the same sign $+$. Namely, these tangles are all PA-tangles. $($If a tangle is alternating and every crossing point in this diagram has the same sign $+$, we call this tangle a PA-tangle.$)$ In addition, the depicted symbols of the tangles as shown in $(2)$ and $(3)$ are such as shown in $(6)$ and those of the tangles as shown in $(4)$ and $(5)$ are such as shown in $(7)$.

\begin{figure}[htbp]
 \centering
 \includegraphics[width=7cm,clip]{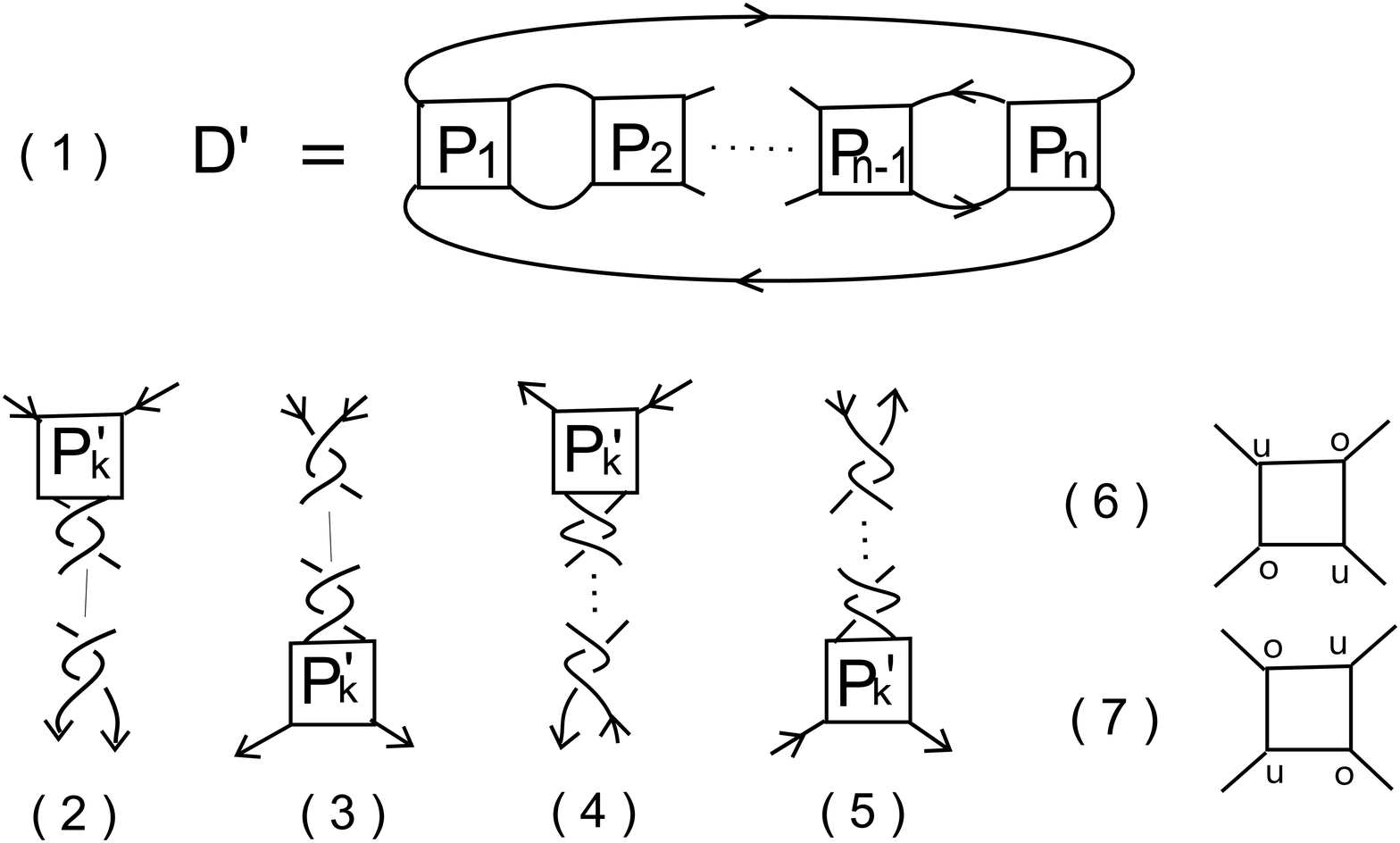}
 \caption{ }
 \label{fig:4-1-10ji}
\end{figure}

Assume $P_s$ is the rightmost tangle of type {\rm II} or type ${\rm III}_+$ as shown in Figure \ref{fig:4-1-10ji}$(2)$ or $(3)$ in the diagram $D'$, then all the tangles situated on the right-hand side of $P_s$ are of type ${\rm III}_-$. Thus when we denote the sum of these tangles by $T_0$, $T_0$ must be a PA-tangle and the depicted symbols of $T_0$ is as shown in Figure \ref{fig:4-1-10ji}$(7)$. Thereby we can transform the diagram $D'$ as shown in Figure \ref{fig:4-1-12ji}, and have a diagram $D''$ which is equivalent to $D'$ where $T_1$ is a PA-tangle.

\begin{figure}[htbp]
  \centering
  \includegraphics[width=10cm,clip]{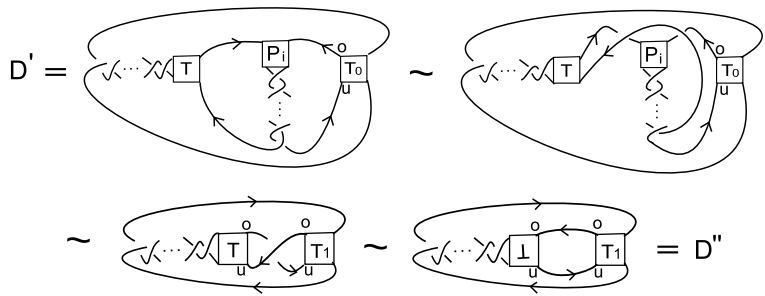}
  \caption{The transformation of diagram D}
  \label{fig:4-1-12ji}
\end{figure}

The directions of outer arcs of $T_1$ are as same as those of $T_0$, and it is obvious that if tangle $P_{s-1}$ which is on the lefthand side of $T_1$ is like as shown in Figure \ref{fig:4-1-10ji}$(4)$ or $(5)$ then the tangle sum $P_{s-1} + T_1$ is a PA-tangle. On the other hand, if $P_{s-1}$ is as in Figure \ref{fig:4-1-10ji}$(2)$ or $(3)$, we can obtain a PA-tangle in a similar fashion like above. Therefore when v is the number of tangles like as shown in Figure \ref{fig:4-1-10ji}$(2)$ or $(3)$, by using v time operations like above we can gain a positive diagram $D'$ as shown in Figure \ref{fig:4-1-14ji} where $T_v$ is a PA-tangle. Moreover we can gain a PA-diagram $D'''$ by changing over/under information of a crossing d. Hence it is clear that $L$ has an almost PA-diagram. This completes the proof of Proposition 4.1.

\qed\end{pr}

\end{document}